\documentclass[10pt]{amsart}

\usepackage{graphicx}

\usepackage[USenglish]{babel}
\usepackage{latexsym}
\usepackage{amsmath}
\usepackage{amsfonts}
\usepackage{amssymb}
\usepackage[all]{xy}
\renewcommand{\a }{\alpha }

\newcommand{\D }{\Delta }

\newcommand{\Sig }{\Sigma}

\newcommand{\be}{\begin{equation}}
\newcommand{\ee}{\end{equation}}
\newenvironment{pf}{\noindent{\sc Proof}.\enspace}{\rule{2mm}{2mm}\medskip}

\newtheorem{remark}{Remark}[section]
\newcommand{\R}{\mathbb{R}}

\newcommand{\Z}{\mathbb{Z}}

\newcommand{\N}{\mathbb{N}}

\newtheorem{theorem}{Theorem}[section]
\newtheorem{proposition}[theorem]{Proposition}
\newtheorem{example}[theorem]{Example}

\newcommand{\bpr}{\begin{proposition}}
\newcommand{\epr}{\end{proposition}}
\newcommand{\bex}{\begin{example}\rm}
\newcommand{\eex}{\end{example}}

\begin{document}

\newtheorem{lem}{Lemma}[section]
\newtheorem{pro}[lem]{Proposition}
\newtheorem{thm}[lem]{Theorem}
\newtheorem{rem}[lem]{Remark}
\newtheorem{cor}[lem]{Corollary}
\newtheorem{df}[lem]{Definition}

\title[On the Leray-Schauder degree of the Toda system]
{On the Leray-Schauder degree of the Toda system on compact
surfaces}

\author{Andrea Malchiodi$^{(1)}$ and David Ruiz$^{(2)}$}

\address{$^{(1)}$ University of Warwick, Mathematics Institute, Zeeman Building, Coventry CV4 7AL and
SISSA, via Bonomea 265, 34136 Trieste (Italy).}

\address{$^{(2)}$
Departamento de An\'alisis Matem\'atico, University of Granada,
18071 Granada (Spain).}

\thanks{A.M. is supported by the FIRB project {Analysis and Beyond},  the PRIN {Variational Methods and
Nonlinear PDE's} and by the University fo Warwick. A. M. and D.R have been supported by the Spanish
Ministry of Science and Innovation under Grant MTM2011-26717. D.
R. has also been supported by J. Andalucia (FQM 116).}

\

\email{A.Malchiodi@warwick.ac.uk, malchiod@sissa.it,
daruiz@ugr.es}

\keywords{Geometric PDEs, Leray-Schauder degree.}

\subjclass[2000]{35J47, 35J61, 58J20.}

\begin{abstract}

In this paper we consider the following {Toda system} of equations
on a compact surface:
$$ \left\{
    \begin{array}{ll}
      - \D u_1 = 2 \rho_1 \left( h_1 e^{u_1}- 1 \right) - \rho_2 \left(h_2 e^{u_2} - 1 \right), \\
     - \D u_2 = 2 \rho_2 \left(h_2 e^{u_2} - 1 \right) - \rho_1 \left(h_1 e^{u_1} - 1 \right). &
    \end{array}
  \right.$$
Here $h_1, h_2$ are smooth positive functions and $\rho_1, \rho_2$
two positive parameters.

In this note we compute the Leray-Schauder degree mod $\Z_2$ of
the problem for $\rho_i \in (4 \pi k, 4 \pi (k+1))$ ($k\in \N$).
Our main tool is a theorem of Krasnoselskii and Zabreiko on the
degree of maps symmetric with respect to a subspace. This result
yields new existence results as well as a new proof of previous
results in literature.

\end{abstract}

\maketitle

\vspace{-0.5cm}

\section{Introduction}

\noindent In this paper we consider the following problem on a
compact orientable surface $\Sigma$:
\begin{equation}\label{eq:e-1}
  \left\{
    \begin{array}{ll}
      - \D u_1 = 2 \rho_1 \left( h_1 e^{u_1}- 1 \right) - \rho_2 \left(h_2 e^{u_2} - 1 \right), \\
     - \D u_2 = 2 \rho_2 \left(h_2 e^{u_2} - 1 \right) - \rho_1 \left(h_1 e^{u_1} - 1 \right). &
    \end{array}
    \right. \end{equation}
Here $h_1, h_2$ are smooth positive functions and $\D$ is the
Laplace-Beltrami operator.

Equation \eqref{eq:e-1} is known as the {Toda system}, and has
been extensively studied in the literature. This problem has a
close relationship to geometry, since it describes the
integrability of Frenet frames for holomorphic curves in
$\mathbb{CP}^2$ (see \cite{gue}). Moreover, it arises in the study
of  non-abelian Chern-Simons theory in the self-dual case, when a
scalar Higgs field is coupled to a gauge potential, see
\cite{dunne, tar, yys}.

Let us first discuss the scalar  counterpart of
\eqref{eq:e-1}, namely a Liouville equation in the form:
\begin{equation} \label{scalar}
 - \Delta u =  \rho \left(h \,  e^{u} - 1 \right),\end{equation} where $\rho \in \R$ and
 $h(x)>0$. Equation \eqref{scalar} arises in the prescribed Gaussian curvature problem
under conformal deformation of the metric, and also describes the
abelian counterpart of \eqref{eq:e-1} from the physical point of
view. This equation has been very much studied in the literature;
there are by now many results regarding existence, compactness of
solutions, bubbling behavior, etc. We refer the interested reader
to the reviews \cite{mreview, tar3}.

Problem \eqref{scalar}  presents a lack of compactness, as its
solutions might blow-up. Indeed, take  a blowing up sequence $u_n$
of \eqref{eq:e-1} with $\rho_n \in \R$ bounded. Then it was proved
in \cite{breme, liyy, ls} that, up to a subsequence
$$ \rho_n \to 8 k \pi, \ k \in \N.$$
Moreover, $e^{u_n}$ behaves like the conformal factor of the
stereographic projection from $S^2$ onto $\R^2$, composed with a
dilation, and located at a finite number of points.

With that result at hand, one can define the Leray-Schauder degree
associated to problem \eqref{eq:e-1} and $\rho \in (8 k \pi, 8
(k+1)\pi)$. By the homotopy property of the degree, the latter is
independent of the metric $g$ and the function $h$, and will
 only vary with $k$ and the topology of $\Sigma$. The computation
of the degree has been accomplished in \cite{clin}, where the
following formula is given:
\begin{equation} \label{grado} d_{LS}= \frac{1}{k!} (-\chi(\Sig)+1) \cdots
(-\chi(\Sig)+k), \ \ (\chi(\Sig)\mbox{ is the Euler characteristic
of } \Sig \mbox{)}.\end{equation} In order to obtain this formula,
in \cite{clin} a detailed study of all blowing-up solutions and
their local degree is performed (see also \cite{mal} for a different approach).

Coming back to  system \eqref{eq:e-1}, it was proved in
\cite{jlw, jw2} that the set of solutions is compact for $
(\rho_1,\rho_2) \notin (4\pi \N \times \R) \cup (\R \times 4\pi
\N)$. In other words, if blowing up occurs, at least one of the
components is quantized.

Therefore, as above, the degree for \eqref{eq:e-1} is well defined
for $(\rho_1, \rho_2)$ away from that set. It is easy to observe
that this degree is equal to $1$ if both $\rho_i$ are smaller than
$4\pi$ (one can deform the parameters to $\rho_1=\rho_2=0$). Apart from
that, there exists no formula for the Leray-Schauder degree for
 system \eqref{eq:e-1} yet.

Because of that, most of the existence results for problem
\eqref{eq:e-1} have used variational methods so far. Indeed, it has been
proved that there exists at least one solution in the following
cases:

\begin{enumerate}

\item for both $\rho_i<4\pi$ (see \cite{jw});

\item for any $\rho_1 < 4\pi$, $\rho_2 \in (4k\pi, 4 (k+1)\pi)$,
$k\in \N$ (see \cite{cheikh});

\item for both $\rho_i \in (4\pi, 8 \pi)$ (see \cite{mr2});

\item for $\rho_1 \in (4k\pi, 4 (k+1)\pi), \rho_2 \in (4m\pi, 4
(m+1)\pi)$, $k,\ m \in \N$ and $\Sigma$ with positive genus (see
\cite{mr3}).

\end{enumerate}

In (1), it is proved that the associated energy functional is
coercive and hence a minimum is found. The rest of the results use
min-max theory, as the functional is no longer bounded from below.

In this note we discuss the parity of the degree
for $\rho_i \in (4 n\pi, 4 (n+1)\pi)$, see Proposition
\ref{p:deg}. Our result is consequence
of a general theorem (recalled in the next section)
concerning the degree of maps symmetric with
respect to a subspace, see \cite{kzabreiko}. We will show that the
degree of the Toda system has the same parity as the degree of the
scalar case with $\rho=\rho_i$, which is given by \eqref{grado}.

In particular, the degree is always odd for $\rho_i \in (4n \pi,
4(n+1)\pi)$ if $n=1,\ 2,\ 3$. The case $n=1$ implies a new,
simpler proof, of the existence result of \cite{mr2}. The cases
$n=2$ or $3$ yield a new existence result:

\begin{theorem}\label{t:main} Assume that $\Sigma$ is homeomorphic to
$\mathbb{S}^2$, and $\rho_i \in (4n \pi, 4(n+1)\pi)$, with $n=2$
or $3$. Then there exists a solution to \eqref{eq:e-1}.
\end{theorem}

\noindent {\bf Acknowledgment:} D. R. thanks Rafael Ortega, from
the University of Granada, for several discussions on the degree
for symmetric maps and his kind help in finding reference
\cite{kzabreiko}.

\section{The parity of the Leray-Schauder degree}

\noindent The main abstract tool we are going to use is the
following one, which is a version of Theorem 21.12 of
\cite{kzabreiko} (page 115).

\begin{thm}\label{t:kz} Let $P$  be a continuous linear projection from a Banach space $E$
onto a (closed) subspace $E^0 \subseteq E$, and define $U=-Id + 2
P$ the reflection with respect to $E^0$, which is assumed to be an
isometry. Let $A : E \to E$ be a compact operator equivariant with
respect to $U$, that is
$$ A \, U(x)= U A(x) \ \ \forall x \in E.$$
Observe that in particular $A(E^0) \subset E^0$. Finally, assume
that $\Phi x = x - A x$ does not vanish on the boundary of
$B_R=B(0,R)$. Then, {\rm deg}$(\Phi,B_R,0)$ and {\rm
deg}$(\Phi|_{E^0}, B_R \cap E^0,0)$ have the same parity.
\end{thm}

\

\begin{remark} Theorem \ref{t:kz} is easy to understand if we
assume that $\Phi$ is $C^1$ and that all its zeroes are
non-degenerate. In such case, all zeroes have  index $\pm 1$,
and the total degree is the sum of the indexes of all zeroes.
Observe now that if $x$ is a zero of $\Phi$, then also $Ux$ is a
zero. Moreover $x=Ux$ if and only if $x \in E^0$. In other words, the
zeroes outside $E^0$ come in pairs, and give an even contribution
to the total degree.

Furthermore, the index of $x \in E^0$ a zero of $\Phi$ could be
different from its index as a zero of $\Phi|_{E^0}$, but in both
cases it is $\pm 1$. So the difference is even.

The proof of \cite{kzabreiko} is topological and does not use
these arguments.

\end{remark}

Let us assume, for simplicity, that $Vol_g(\Sigma) = 1$. For $\a
\in (0,1)$ we consider the functional framework
$$
 E = C^{2,\alpha}_0(\Sigma) \times  C^{2,\alpha}_0(\Sigma),
$$
where $C^{2,\alpha}_0(\Sigma)$ stands for the class of $C^{2,\a}$
functions with zero average. Define now the operator
$A = A_{\rho_1,\rho_2}^{h_1,h_2} : E \to E$ as
$$
A \left( \begin{array}{c}
  u_1 \\ u_2
  \end{array} \right) =
\left( \begin{array}{c}
  (-\D)^{-1} \left[ 2 \rho_1 \left( h_1 \frac{e^{u_1}}{\int_{\Sigma} h_1 e^{u_1} dV_g}
  - 1 \right) - \rho_2 \left(h_2 \frac{e^{u_2}}{\int_{\Sigma} h_2e^{u_2} dV_g}- 1 \right) \right]
  \\ (-\D)^{-1} \left[ 2 \rho_2 \left(h_2 \frac{e^{u_2}}{\int_{\Sigma} h_2 e^{u_2} dV_g} - 1 \right)
  - \rho_1 \left(h_1 \frac{e^{u_1}}{\int_{\Sigma} h_1e^{u_1} dV_g}  - 1 \right) \right]
  \end{array} \right).
$$
Here by $(-\D)^{-1} f$, $f \in C^\a(\Sigma)$, we denote the unique
solution $u$ of $- \D u = f$ with zero average. In the above
formula solutions exist by Fredholm's theory. Notice also that
zeroes of $ \Phi= Id - A$ give rise to solutions of
\eqref{eq:e-1}. Indeed, it suffices to add proper constants to
$u_1, u_2$ in order to  have $\int_{\Sigma} h_i e^{u_i} dV_g = 1$.

\

By elliptic regularity theory the operator $A$ is compact.
Moreover, if $n \in \N$ and $\rho_1$, $\rho_2 \in (4 n \pi,
4(n+1)\pi)$, the solutions are a priori bounded, see \cite{jlw}.
Therefore for such values of $\rho_1, \rho_2$ and for $R$
sufficiently large the degree {\rm deg}$(\Phi,B_R,0)$ is well
defined. The main result of this paper is the following:

\begin{pro}\label{p:deg}
Consider $n\in \N$ and let $\rho_1$, $\rho_2 \in (4n\pi, 4(n+1)\pi)$.
Then, for sufficiently large $R$, {\rm deg}$(\Phi,B_R,0)$ has the
same parity as
$$d_k= \frac{1}{k!} (-\chi(\Sig)+1) \cdots
(-\chi(\Sig)+k),$$ where $\chi(\Sig)$ is the Euler characteristic
of $\Sig$ and $k= [n/2]$.
\end{pro}

\begin{pf}
Let us consider the homotopy
$$
  (1-s) (\rho_1,\rho_2) + s (\rho, \rho); \qquad \quad \rho = \frac{1}{2} (\rho_1 + \rho_2);
$$
$$
  (1-s) (h_1,h_2) + s (h, h); \qquad \quad h = \frac{1}{2} (h_1 +
  h_2).
$$
By the homotopy invariance of the degree we deduce that (for $R$
sufficiently large)
\begin{equation}\label{eq:degsame1}
  {\rm deg}(\Phi_{\rho_1,\rho_2}^{h_1,h_2},B_R,0) = {\rm deg}(\Phi_{\rho,\rho}^{h,h},B_R,0).
\end{equation}
Because of that, it suffices to study the degree for $h_1=h_2=h$,
$\rho_1=\rho_2= \rho$.

We choose $E^0$ to be the couples of identical functions in $E$,
namely
$$
  E^0 = \left\{ \left( \begin{array}{c}
  u_1 \\ u_2
  \end{array} \right) \in E \; : \; u_1 = u_2 \right\},
$$
and we define the projection $P : E \to E^0$ as
$$
  P \left( \begin{array}{c}
    u_1 \\ u_2
    \end{array} \right) = \frac{1}{2}  \left( \begin{array}{c}
      u_1+u_2 \\ u_1+u_2
      \end{array} \right).
$$
Observe that the reflection $U$ is given by:
$$
  U \left( \begin{array}{c}
    u_1 \\ u_2
    \end{array} \right) =   \left( \begin{array}{c}
      u_2 \\ u_1
      \end{array} \right).
$$

With these definitions, we are in conditions to apply Theorem
\ref{t:kz}, and hence ${\rm deg}(\Phi, B_R,0)$ has the same parity
as ${\rm deg}(\Phi|_{E_0},B_R\cap E_0,0)$.

We now define:

$$T: C^{2,\alpha}_0(\Sigma) \to E_0,\ T(u)=\left( \begin{array}{c}
    u \\ u
    \end{array} \right) .$$
Clearly $T$ is an homeomorphism, which implies that:
$$ {\rm deg}(\Phi|_{E_0},B_R\cap E_0,0)= {\rm deg}(T^{-1} \circ \Phi|_{E_0} \circ T, \tilde{B}_R,0),$$ where $\tilde{B}_R=B(0,R) \subset C^{2,\alpha}_0(\Sigma)$.

Observe now that:

$$\tilde{\Phi}(u):=T^{-1} \circ \Phi|_{E_0} \circ T(u)= u- (-\D)^{-1} \left[  \rho \left( h \frac{e^{u}}{\int_{\Sigma} h e^{u} dV_g}
  - 1 \right)\right].$$
Moreover, $\rho \in (4n\pi, 4(n+1)\pi) \subset (8k\pi, 8(k+1)\pi)$
for $k=[n/2]$. Finally, ${\rm deg}(\tilde{\Phi}, \tilde{B}_R,0)$
has been computed in \cite{clin} and it is given by the formula
\eqref{grado}. Notice that the Leray-Schauder degree in
\cite{clin} has been computed in the $H^1$ setting, but using
elliptic regularity theory one can prove that this coincides with
the degree in the $C^{2,\a}$ setting, see Theorem B.1 in
\cite{yy95}. This concludes the proof.
\end{pf}

\

\noindent Observe that if $n=1$, then $k=0$ and hence the degree
is odd. In this way we obtain the existence result of \cite{mr2}
with an alternative approach. There are other cases in which the
degree is odd, so we recover some of the results of \cite{mr3}. In
particular, if $n=2,3$, then $k=1$ and the degree is odd for any
compact orientable $\Sigma$. As a consequence, we obtain Theorem
\ref{t:main}, which gives a new existence result.


\begin{thebibliography}{99}

\bibitem{mr3}{Luca Battaglia, Aleks Jevnikar, Andrea Malchiodi and David Ruiz,}{ A general
existence result for the toda system on compact surfaces, }{
preprint arXiv 1306.5404, 2013.}

\bibitem{breme}
Haim Brezis and Frank Merle.
\newblock Uniform estimates and blow-up behavior for solutions of {$-\Delta
  u=V(x)e^u$} in two dimensions.
\newblock {Comm. Partial Differential Equations}, 16(8-9):1223--1253, 1991.



\bibitem{clin}{Chiun-Chuan Chen and Chang-Shou Lin, }{Topological degree for a mean field
equation on Riemann surfaces, }{Comm. Pure Appl. Math.,
56(12):1667-1727, 2003.}



\bibitem{dunne}
G. Dunne.
\newblock {Self-dual Chern-Simons Theories}.
\newblock Lecture notes in physics. New series m: Monographs. Springer, 1995.

\bibitem{gue}{M. A. Guest, }{Harmonic maps, loops groups, and integrable systems. }{London Mathematical
Society Student Texts, 38. Cambridge University Press, Cambridge,
1997.}

\bibitem{jlw}
J{\"u}rgen Jost, Changshou Lin, and Guofang Wang.
\newblock Analytic aspects of the {T}oda system. {II}. {B}ubbling behavior and
  existence of solutions.
\newblock {Comm. Pure Appl. Math.}, 59(4):526--558, 2006.

\bibitem{jw}
J{\"u}rgen Jost and Guofang Wang.
\newblock Analytic aspects of the {T}oda system. {I}. {A} {M}oser-{T}rudinger
  inequality.
\newblock {Comm. Pure Appl. Math.}, 54(11):1289--1319, 2001.

\bibitem{jw2}
J{\"u}rgen Jost and Guofang Wang.
\newblock Classification of solutions of a {T}oda system in {${\mathbb R}^2$}.
\newblock {Int. Math. Res. Not.}, (6):277--290, 2002.

\bibitem{kzabreiko}{M. A. Krasnoselskii and P. P. Zabreiko, }{Geometrical methods of nonlinear anal-
ysis, }{volume 263 Fundamental Principles of Mathematical
Sciences. Springer-Verlag, Berlin, 1984. Translated from the
Russian by Christian C. Fenske.}

\bibitem{yy95}{Yan Yan Li, }{Prescribing scalar curvature on $S^n$ and related problems. I.,}{J.
Differential Equations, 120(2):319-410, 1995.}

\bibitem{liyy}
Yan~Yan Li.
\newblock Harnack type inequality: the method of moving planes.
\newblock {Comm. Math. Phys.}, 200(2):421--444, 1999.

\bibitem{ls}
Yan~Yan Li and Itai Shafrir.
\newblock Blow-up analysis for solutions of {$-\Delta u=Ve^u$} in dimension
  two.
\newblock {Indiana Univ. Math. J.}, 43(4):1255--1270, 1994.

\bibitem{mal}
Andrea Malchiodi.
\newblock Morse theory and a scalar field equation on compact surfaces.
\newblock {Adv. Differential Equations}, 13(11-12):1109--1129, 2008.

\bibitem{mreview}
Andrea Malchiodi.
\newblock Topological methods for an elliptic equation with exponential
  nonlinearities.
\newblock {Discrete Contin. Dyn. Syst.}, 21(1):277--294, 2008.

\bibitem{cheikh}
Andrea Malchiodi and Cheikh~Birahim Ndiaye.
\newblock Some existence results for the {T}oda system on closed surfaces.
\newblock {Atti Accad. Naz. Lincei Cl. Sci. Fis. Mat. Natur. Rend. Lincei
  (9) Mat. Appl.}, 18(4):391--412, 2007.

\bibitem{mr2}
Andrea Malchiodi and David Ruiz.
\newblock A variational analysis of the {T}oda system on compact surfaces.
\newblock {Comm. Pure Appl. Math.}, 66(3):332--371, 2013.


\bibitem{tar}
Gabriella Tarantello.
\newblock {Selfdual gauge field vortices}.
\newblock Progress in Nonlinear Differential Equations and their Applications,
  72. Birkh\"auser Boston Inc., Boston, MA, 2008.
\newblock An analytical approach.

\bibitem{tar3}
Gabriella Tarantello.
\newblock Analytical, geometrical and topological aspects of a class of mean
  field equations on surfaces.
\newblock {Discrete Contin. Dyn. Syst.}, 28(3):931--973, 2010.



\bibitem{yys}
Yisong Yang.
\newblock {Solitons in field theory and nonlinear analysis}.
\newblock Springer Monographs in Mathematics. Springer-Verlag, New York, 2001.



\end{thebibliography}
\end{document}